\documentclass{amsart}
\usepackage[english]{babel}
\usepackage{latexsym}
\usepackage{amssymb}
\usepackage{amsmath}
\DeclareMathOperator{\Ker}{Ker} 
\DeclareMathOperator{\TrueIm}{Im} \DeclareMathOperator{\ord}{ord}
\DeclareMathOperator{\Aut}{Aut} \DeclareMathOperator{\Int}{Int}
\DeclareMathOperator{\Out}{Out} \DeclareMathOperator{\Lie}{Lie}
\DeclareMathOperator{\idd}{id} \DeclareMathOperator{\diagg}{diag}
\newtheorem{theor}{$\phantom{sss}$Theorem}
\newtheorem{utver}{$\phantom{sss}$Proposition}
\newtheorem{laemma}{$\phantom{sss}$Lemma}
\newtheorem{corol}{$\phantom{sss}$Corollary}
\newtheorem{fiend}{$\phantom{sss}$Remark}
\newtheorem{examp}{$\phantom{sss}$Example}

\begin{document}
\fontfamily{ptm} \fontsize{12pt}{16pt} \selectfont
\author[S.N. Fedotov]{Stanislav N. Fedotov}
\title[Groups with periodic components]{Affine algebraic groups with periodic components}
\thanks{Supported by grant NSh-1983.2008.1}
\address{Moscow State University, department of Higher Algebra
\endgraf Email: glwrath@yandex.ru}

\maketitle
\begin{abstract} A connected component of an affine
algebraic group is called periodic if all its elements have finite
order. We give a characterization of periodic components in terms of
automorphisms with finite number of fixed points. It is also
discussed which connected groups have finite extensions with
periodic components. The results are applied to the study of the
normalizer of a maximal torus in a simple algebraic group.
\end{abstract}

\section{Introduction}

It is well known that every connected affine complex algebraic group
of positive dimension contains elements of infinite order. On the
other hand, for a non-connected group some of connected components
may consist of elements of finite order. Such components we call
{\it periodic}. For example, one of the connected components of the
orthogonal group $O_2(\mathbb{C})$ consists of reflections and
therefore is periodic. This work is devoted to the study of affine
algebraic groups with periodic components.

Let $G$ be an affine algebraic group over an algebraic closed field
$\Bbbk$ of characteristic zero and $G^0 = H$ be its connected
component of unity. In Section 2 we formulate equivalent conditions
for a component $gH$ of the group $G$ to be periodic in terms of
automorphisms induced by the action of the elements of $gH$ on $H$
by conjugation. Namely, a component $gH$ is periodic if and only if
the automorphism $\varphi_g: H\rightarrow H,\ h\mapsto g^{-1}hg$ has
only finitely many fixed points. Also we prove that for a periodic
component $gH$ the action of $H$ on $gH$ by conjugation is
transitive and, consequently, all elements of $gH$ have same orders.
In Section 3 we investigate which connected algebraic groups have
finite extensions with periodic components and show that all such
groups are solvable. It turns out that solvability is only a
necessary condition: we give an example of a series of connected
solvable groups having no extensions with periodic components.
Section 4 is devoted to the study of torus extensions and, in
particular, to estimates of order of elements in periodic components
of such extensions. Namely, if there is an extension $G = T\cup
gT\cup\ldots\cup g^{m-1}T$ of a torus $T$ with a periodic component
$gT$, and $\ord(\varphi_g) = k$, then $\ord(g)$ divides $mk$. using
the obtained results in Section 5 we investigate periodic components
of the normalizer of a maximal torus in simple groups and find the
order of elements in these components. In particular, we give a
formula for the number of periodic components in the normalizer and
show that the components defined by the Coxeter elements of the Weyl
group are periodic.

The author thanks I.V. Arzhantsev for the idea of the work and
useful discussions, D.A. Shmelkin for a preliminary version of
Theorem 1 and E.I. Khukhro for bringing the paper [7] to his
attention.

\section{Periodic components and automorphisms}

The following theorem suggests several characterizations of periodic
components.

\begin{theor} Let $G$ be an affine algebraic group and $H = G^0$ be
its connected component of unity. For each $g\in G$ define the
automorphism $\varphi_g : H\longrightarrow H; h\mapsto g^{-1}hg$.
Then the following conditions are equivalent:

(1) the connected component $gH$ is periodic;

(2) the subgroup of fixed points of the automorphism $\varphi_g$ is
finite;

(3) the action of $H$ on $gH$ by conjugation is transitive.
\end{theor}

\par

\begin{proof} First of all note that $H^{\varphi_g}$ as a subgroup
of fixed points of a regular automorphism is an algebraic subgroup.

\medskip

{\it $(1)\Rightarrow (2)$} Let the component $gH$ consist of
elements of finite order with the group $H^{\varphi_g}$ being
infinite. Then the latter has positive dimension and therefore
contains an algebraic subgroup of dimension $1$. But such a subgroup
may only be either a one-dimensional torus or an additive group of
the ground field \cite[Th.~3.2.8]{Vind} and therefore contains an
element of infinite order. Denote it by $h_0$.

\medskip

But $h_0g = gh_0$ and so $\forall n\in\mathbb{N}: (gh_0)^n =
g^nh_0^n\ne e$, since $g^n$, as an element of finite order, can not
be inverse to $h_0^n$. This yields a contradiction.

\medskip

{\it $(2)\Rightarrow (3)$} Let $\left|H^{\varphi_g}\right| <\infty$.
The quotient group $G/H$ is finite and, consequently, $\exists
m\in\mathbb{N}$ such, that $g^m\in H$. But then $g^m\in
H^{\varphi_g}$, implying that $g$ has finite order.

\medskip

Consider the morphism $f_g: H\longrightarrow gH, h\mapsto hgh^{-1}$
(it is well-defined, because $\forall h\in H: hgh^{-1} =
g\varphi_g(h)h^{-1}\in gH$).

\medskip

\begin{laemma} If the subgroup $H^{\varphi_g}$ is finite, then the fibers
of the morphism $f_g$ are finite.
\end{laemma}

\begin{proof} Fix an arbitrary
$h_0\in H$. Then $\forall h\in H: hgh^{-1} =
h_0gh_0^{-1}\Leftrightarrow h_0^{-1}hgh^{-1}h_0 = g\Leftrightarrow
h_0^{-1}h = gh_0^{-1}hg^{-1}\Leftrightarrow h_0^{-1}h\in
H^{\varphi_g}$, implying $h\in h_0H^{\varphi_g}$. But the subgroup
$H^{\varphi_g}$ is finite.

\end{proof}

\medskip

Denote by $C(g)$ the image of the morphism $f_g$. By the theorem on
the dimension of fibers of a morphism  $\dim{C(g)} = \dim{H}$
\cite[Ch.~I \S~6.3]{Antisem}. But, as a connected component of an
algebraic group, $gH$ is irreducible, and therefore $C(g)$ is dense
in $gH$. This yields that the component $gH$ contains dense subset,
consisting of the elements conjugate to $g$ and, consequently, of
the same order.

\medskip

Let $N$ be the order of the element $g$. Note that the subset $M =
\{x\in gH\,|\,(gh)^N = e\}$ is closed in $gH$ in Zariski topology.
On the other hand, $M$ contains the subset $C(g)$ that is dense in
$gH$, and therefore $M$ is equal to $gH$.

So $\forall s\in gH$ the order of $s$ is finite and, consequently,
the image $C(s)$ of the morphism $f_s$ is dense in $gH$ and is
opened in $gH$ as an orbit of $H : gH$. But for different $s\in gH$
the sets $C(s)$ either are disjoint or are equal, and as opened
subsets, they must have common points. This means that $gH =
C(s)\,\forall s\in gH$.

{\it $(3)\Rightarrow (1)$} Consider the Jordan decomposition $g =
tu$ with $t$ being semi-simple, $u$ being unipotent and $tu = ut$.
Then $t,u\in G$ \cite[Th.~3.4.6]{Vind} and $G(u)\subset G$ where
$G(u)$ is the minimal algebraic subgroups containing $u$. But $G(u)$
is unipotent and therefore connected implying $G(u)\subset H$, so
$u\in H$. Consequently, $t = gu^{-1}\in gH$ and, as all elements in
$gH$ are conjugate, $g$ is also semi-simple.

Assume that $g$ has infinite order. Then the component $gH$ contains
infinitely many elements that are its powers. Consider the faithful
linear representation $\rho : G\longrightarrow GL(V)$
~\cite[Th.~3.1.8]{Vind}. As $g$ is semi-simple, it and all its
powers are represented by diagonal matrices in an appropriate basis
of $V$. But those of them lying in $\rho(gH)$ are conjugate (by
matrices of $\rho(H)$) and consequently may differ only by a
permutation of diagonal elements. So among the degrees of the matrix
$\rho(g)$ that lye in $\rho(gH)$ there are only finitely many
different ones. But the representation $\rho$ is faithful, thus
yielding the contradiction.

\end{proof}

\begin{corol} If the component $gH$ is periodic, then all its
elements have same orders.
\end{corol}

\begin{corol} If the component $gH$ is periodic, then $H$ lyes
in the commutant $[G, G]$ of the group $G$.
\end{corol}

\begin{proof} The condition
$gH = Hg = \{hgh^{-1}\mid h\in H\}$ implies $H =
\{hgh^{-1}g^{-1}\mid h\in H\}\subseteq [G, G]$.

\end{proof}

So if the group $G$ contains a periodic component, then computation
of its commutant can be reduced to computation of the commutant of
the finite group of its components.

\begin{corol} If the group $G$ has a periodic component, then
for each its linear representation $\rho: G\longrightarrow GL(V)$
the image of $\rho(G^0)$ lyes in $SL(V)$.
\end{corol}

Now let $H$ be a connected group and $\varphi$ be its automorphism
of order $m < \infty$ with finitely many fixed points. Consider the
semi-direct product $H\,\leftthreetimes\langle a\rangle_m$ defined
by this automorphism, where $\langle a\rangle_m$ is a cyclic group
of order $m$ with generator $a$. By Theorem 1 the coset $\{(h,
a)\mid h\in H\}$ is a periodic component of the resulting group.
Thus we get the following criterion: a connected group has a finite
extension with a periodic component if and only if it has an
automorphism of finite order with finitely many fixed points.

Let $G$ be a finite extension of a connected group $H$. Note that
$\forall g\in G$ the automorphism $\varphi_g$ has finite order in a
quotient group $\Aut(H)/\Int(H)$. In fact, in the previously
described criterion we can replace the condition for an automorphism
to have finite order by the condition to have finite order modulo
the group $\Int(H)$.

\begin{utver} Let $H$ be a connected algebraic group, $\varphi$ be its automorphism,
of order $k < \infty$ modulo inner automorphisms of $H$ and such,
that $|H^{\varphi}| <\infty$. Then the order of automorphism
$\varphi$ is finite.
\end{utver}
\begin{proof} We have $\varphi^k\in \Int{H}$ and,
consequently, $\exists x\in H$ such, that $\varphi^k(h) = xhx^{-1}\
\forall h\in H$. For an arbitrary $h\in H$ compute $\varphi^{k +
1}(h)$ in two ways: on one hand
$$\varphi^{k + 1}(h) = \varphi^k(\varphi(h)) = x\varphi(h)x^{-1},$$
and on the other hand
$$\varphi^{k + 1}(h) = \varphi(\varphi^k(h)) = \varphi(xhx^{-1}) =
\varphi(x)\varphi(h)(\varphi(x))^{-1}.$$

But $\varphi$ is an automorphism, therefore $\varphi(H) = H$ and
$\forall h\in H$ we have $x\varphi(h)x^{-1} =
\varphi(x)\varphi(h)\varphi(x)^{-1}$, which means that
$x^{-1}\varphi(x)$ lyes in the center $Z(H)$ of the group $H$. So we
get $x^{-1}\varphi(x) = xx^{-1}\varphi(x)x^{-1} = \varphi(x)x^{-1}$.
Consequently, $\forall s\in\mathbb{N}:$
\begin{equation*}
x^{-s}\varphi(x^s) = x^{-s}(\varphi(x))^s = (x^{-1}\varphi(x))^s\in
Z(G).
\end{equation*}

Consider the morphism $f: H\rightarrow H, h\mapsto
h^{-1}\varphi(h)$.

\begin{laemma} Fibers of the morphism $f$ are finite.
\end{laemma}

\begin{proof} Fix an element
$h_0\in H$. Then $\forall g\in H$ one has: $h_0^{-1}\varphi(h_0) =
g^{-1}\varphi(g)\Leftrightarrow gh_0^{-1} =
\varphi(g)\varphi(h_0^{-1}) = \varphi(gh_0^{-1})\Leftrightarrow
gh_0^{-1}\in H^{\varphi}\Leftrightarrow g\in h_0H^{\varphi}$. But
the group $H^{\varphi}$ is finite.
\end{proof}

Consider morphisms $f_s: x^sZ(H)\rightarrow Z(H), h\mapsto
h^{-1}\varphi(h), s = 1,2,...$. As $f_s = f|_{x^sZ(H)}$ fibers of
$f_s$ are finite and contain less than $|H^{\varphi}|$ points. On
the other hand by the theorem on the dimension of fibers of a
morphism $\dim{\TrueIm f_s} = \dim{x^sZ(H)} = \dim{Z(H)}$. Hence for
all $s$ the set $\overline{\TrueIm f_s}$ is a union of several
connected components of the group $Z(H)$. Therefore, by the theorem
on the image of a dominant morphism $\TrueIm f_s$ contains an open
in Zariski topology subset of this union. If the number of different
cosets $x^sZ(H)$ is sufficiently big, then we get the contradiction
with a limitation on the number of points in a fiber of the morphism
$f$. Namely, let $d$ be a number of connected components of $Z(H)$,
$m$ be less than $d\cdot|H^{\varphi}|$ and all $x^sZ(H), s =
1,2,...,m$ be different. Then among the sets $\TrueIm{f_s}\ (s =
1,2,...m)$ some $|H^{\varphi}| + 1$ have non-empty intersection. Let
$g_0$ be an element lying in $|H^{\varphi}| + 1$ of the sets. Then
$f^{-1}(g_0) \supseteq \cup_{s = 1}^{m}f_s^{-1}(g_0)$ and
consequently $|f^{-1}(g_0)| \geqslant m/d
> |H^{\varphi}|$, which is a contradiction.

Hence among the cosets $x^sZ(H), s = 1, 2,...\,$ only finitely many
of them are different and therefore $\exists q:\,x^q\in Z(H)$ (it
was shown that $q\leqslant d\cdot |H^{\varphi}|$ where $d$ is a
number of connected components of $Z(H)$). But then $\varphi^{kq} =
\idd$.

\end{proof}

\section{Extensions with periodic components}

Now we turn to the question, whether a given connected algebraic
group has extensions with periodic components.

\begin{utver} If a reductive algebraic group has a periodic
component, then its connected component of unity is a torus.
\end{utver}

\begin{proof} Let $G$ be a reductive group and $gG^0$  be its periodic component.
Since every connected reductive group $G^0$ is an almost direct
product of a central torus $T$ and a semi-simple subgroup
$S$\cite[Ch.~6]{Vind}. But $S = [G^0, G^0]$ and consequently it is
stable under the action of the automorphism $\varphi_g$. Hence it is
sufficient to prove the following lemma.

\begin{laemma} The subgroup of fixed points of an automorphism of a
semi-simple group is infinite.
\end{laemma}

\begin{proof} Let $H$ be a connected semi-simple group and $\mathfrak{h} =
\Lie{H}$ be its Lie algebra. Consider a map $D:
\Aut{H}\longrightarrow \Aut{\mathfrak{h}},\; f\mapsto d_ef$. As
$D(\Int{H}) = \Int{\mathfrak{h}}$, it is evident that $D$ defines
correctly a mapping from $\Out{H} = \Aut{H}/\Int{H}$ to
$\Out{\mathfrak{h}} = \Aut{\mathfrak{h}}/\Int{\mathfrak{h}}$.
For semi-simple Lie algebras the following fact is well known.

\begin{utver} The group $\Int{\mathfrak{h}}$ is a connected component of
unity in $\Aut{\mathfrak{h}}$; diverse connected components of the
group $\Aut{\mathfrak{h}}$ are the sets
$(\Int{\mathfrak{h}})\hat\tau$ for diverse $\tau\in \Aut{\Pi}$
($\Pi$ is a positive roots system of algebra $\mathfrak{h}$), where
$\hat\tau$ are defined in the following way:
$$\hat\tau(h_{\alpha}) = h_{\tau^{-1}(\alpha)},\;
\hat\tau(e_{\alpha}) = e_{\tau^{-1}(\alpha)},\; \hat\tau(e_{ -
\alpha}) = e_{ - \tau^{-1}(\alpha)}\quad (\alpha\in\Pi).$$
\end{utver}

(Here $h_{\alpha}, e_{\alpha}, e_{ - \alpha}$ are the canonic
generators of the algebra $\mathfrak{h}$; for more information
see~\cite[\S~4, Ch.~4]{Vind}.)

\medskip

Hence for an appropriate choice of a maximal torus and a positive
system $\Pi$ the automorphism $d_e\varphi$ acts as $\hat\tau$ for a
certain $\tau\in\Aut\Pi$. But then consider $\lambda =
\sum\limits_{\alpha\in\Pi} h_{\alpha}$. By definition of $\hat\tau$
it is evident that $d_e\varphi$ acts identically on the line spanned
by $\lambda$ and therefore $\varphi$ acts identically on the
corresponding one-dimensional subgroup. This implies
$\left|H^{\varphi}\right| = \infty$.
\end{proof}

This proves Proposition 2.

\end{proof}

\begin{theor} If an affine algebraic group $G$ has a periodic component,
then $G^0$ is solvable.
\end{theor}

\begin{proof} Consider the Levi decomposition: $G^0 = L\rightthreetimes U$,
where $U$ is the unipotent radical, and $L$ is reductive (i.e. a
Levi subgroup)\cite[\S~4, Ch.~6]{Vind}. Let $gG^0$ be a periodic
component. Then $gLg^{-1}$ is also a Levi subgroup and by Maltsev's
Theorem $\exists h\in U:\, gLg^{-1} = hLh^{-1}$; hence the
automorphism $\varphi_{h^{-1}g}$ stabilizes the subgroup $L$. But as
$L$ is reductive, Proposition 2 implies that the restriction of the
automorphism $\varphi_{h^{-1}g}$ on $L$ may have finitely many fixed
points if and only if $L$ is a torus. This means that $G^0 =
L\rightthreetimes U$ is solvable.

\end{proof}

The automorphism of a Lie algebra is called {\it regular} provided
it fixes no point except the zero. Note that an automorphism of an
algebraic group having finitely many fixed points induces a regular
automorphism of its Lie algebra. After having proved Theorem 2 the
author found the work~\cite{Borel} where it is proved that a finite
dimensional Lie algebra over the field of characteristic zero
possessing a regular automorphism of finite order is solvable. It
gives an alternative proof for Theorem 2. In ~\cite{Krec} one can
find some other properties of Lie algebras with regular
automorphism.

\begin{examp} {\rm Consider $G = U_n\cup gU_n\cup\ldots\cup g^nU_n$
where $U_n\subset GL_n(\Bbbk)$ as the group of all uni-triangular
matrices and $g$ is a diagonal matrix with eigenvalues $1, \xi,
\xi^2,\ldots, \xi^{n - 1}$, where $\xi$ is a primitive root of unity
of degree $n$. All elements of the connected components $gU_n$ are
upper-triangular matrices with different roots of unity on the
diagonal, hence they are semi-simple and have finite order.}
\end{examp}

Now the problem arises that is if every unipotent group has
extensions with periodic components or, just the same, if every
nilpotent Lie algebra possesses a periodic regular automorphism. By
improved Ado's Theorem~\cite[Ch. 1, \S 5.3]{Vindex} every finite
dimensional nilpotent Lie algebra is isomorphic to a subalgebra of
the algebra of all nilpotent triangular matrices of a certain
dimension $n\in \mathbb{N}$. By the Campbell-Hausdorff
Formula~\cite[Part~1, Ch.~IV, \S~7,8]{Serr} the image of this
subalgebra under the exponential mapping is a subgroup in $U_n$ and,
furthermore, the exponential mapping establishes an isomorphism of
algebraic manifolds. So if there exists a nilpotent Lie algebra
$\mathfrak{g}$ with only unipotent automorphisms, then the
corresponding unipotent group will have no extensions with periodic
components. An example of such an algebra can be found
in~\cite{Dyer}.

The following proposition shows that if both torus $T$ and
unipotent group $U$ have extensions with periodic components, then
their semidirect product may have no such extensions.

\begin{utver} Among the algebraic groups of type $\Bbbk^*\rightthreetimes\Bbbk$
only $\Bbbk^*\times\Bbbk$ has a finite extension with periodic
components.
\end{utver}

\begin{proof} In case if the product is not direct it is sufficient to
note that a non-commutative two-dimensional Lie algebra has no
regular automorphisms of finite order.

For the group $\Bbbk^*\times\Bbbk$ there is an extension with
periodic components. It can be constructed as follows:

\smallskip

$$G = \left\{
\begin{pmatrix}
\begin{matrix}
t & 0\\
0 & t^{-1}
\end{matrix} & \mbox{\huge 0}\\
\mbox{\huge 0} &
\begin{matrix}
1 & s\\
0 & 1
\end{matrix}
\end{pmatrix}\right\}\cup\left\{
\begin{pmatrix}
\begin{matrix}
0 & t^{-1}\\
t & 0
\end{matrix} & \mbox{\huge 0}\\
\mbox{\huge 0} &
\begin{matrix}
1 & s\\
0 & -1
\end{matrix}
\end{pmatrix}\right\}\ ,t\in\Bbbk^*, s\in\Bbbk.
$$

\end{proof}

One may suggest that only the groups of type $T\times U$, where $T$
is a torus and $U$ is a unipotent group, have automorphisms with
finitely many fixed points. But it is not true.

\begin{examp} {\rm Consider the group $H = T\rightthreetimes U$,
where $T\cong(\Bbbk^*)^{n-1}, U\cong\Bbbk^n$, and multiplication is
defined as follows: $(t,u)\cdot(s,v)=(ts,\psi(s)(u)+v)$, where
$$\psi(s)(u)=(s_1^{-1}u_1,\ldots,s_{n-1}^{-1}u_{n-1},s_1\ldots s_{n-1}u_n).$$
Fix a natural number $k$, that divides $n$, and define the mapping
$\varphi: H\rightarrow H$ as follows: $\varphi((t,u))=(\beta(t),
\alpha(u))$, where $\beta((t_1, t_2, \ldots t_{n-1})) = (t_2,
\ldots, t_{n-1}, t_1^{-1}t_2^{-1}\ldots t_{n-1}^{-1})$,
$$\alpha(u)=
\begin{pmatrix}
0 & 1 &&&\\
& 0 & 1 &&\\
&& \ddots & \ddots &\\
&&& 0 & 1\\
\xi &&&& 0
\end{pmatrix}u,$$
where $\xi$ is a primitive root of unity of degree $k$. We should
now prove that $\varphi$ is an automorphism:

$$\varphi((t,u)(s,v))=\varphi((ts,\psi(s)(u)+v))=
(\beta(ts), \alpha(\psi(s)(u))+\alpha(v));$$
$$\varphi((t,u))\varphi((s,v))=(\beta(t),\alpha(u))(\beta(s),
\alpha(v))=(\beta(t)\beta(s),\psi(\beta(s))(\alpha(u))+\alpha(v)).$$
But it is easy to see, that $\forall t\in T, u\in U:$
$$\alpha(\psi(t)(u))=$$
$$ = \begin{pmatrix}
0 & 1 &&&\\
& 0 & 1 &&\\
&& \ddots & \ddots &\\
&&& 0 & 1\\
\xi &&&& 0
\end{pmatrix}
\begin{pmatrix}
t_1^{-1} &&&& 0\\
 & t_2^{-1} &&&\\
 && \ddots &&\\
 &&& t_{n-1}^{-1} &\\
 0 &&&& t_1t_2\ldots t_{n-1}
\end{pmatrix}u =$$
$$=\begin{pmatrix}
 t_2^{-1} &&&& 0\\
 & \ddots &&&\\
 && t_{n-1}^{-1} &&\\
 &&& t_1t_2\ldots t_{n-1} &\\
 0 &&&& t_1^{-1}
\end{pmatrix}
\begin{pmatrix}
0 & 1 &&&\\
& 0 & 1 &&\\
&& \ddots & \ddots &\\
&&& 0 & 1\\
\xi &&&& 0
\end{pmatrix}u= $$
$$= \psi(\beta(t))(\alpha(u)),$$ so $\varphi$ is an automorphism.
Furthermore, $\varphi$ has order $nk$ and as many as $n$ fixed
points $((\eta, \ldots, \eta), 0)$, where $\eta$ is a root of unity
of degree $n$.

We now construct the matrix realization of the group $H$, that
realizes the automorphism $\varphi$:
$$H = \left\{
\begin{pmatrix}
T & K\\
0 & E
\end{pmatrix}\left|
\begin{matrix}
T = \diagg(t_1,\ldots, t_{n-1}, t_n),\ 
t_i\in\Bbbk^*\\
\ K = \diagg(t_1a_1,\ldots, t_{n-1}a_{n-1}, t_na_n),\ a_i\in\Bbbk\\
t_1\ldots t_n = 1
\end{matrix}\right.\right\}$$

$$G = H\cup g_0H\cup\ldots\cup\ g_0^{n - 1}H, \mbox{ где}$$
$$g_0 =
\begin{pmatrix}
\Sigma & 0\\
0 & \Sigma D
\end{pmatrix},\
\Sigma = \begin{pmatrix}
0 &&&&\xi^{n-1}\\
1 & 0 &&&\\
& \ddots & \ddots &&\\
&& 1 & 0 &\\
 &&& 1& 0
\end{pmatrix},\ D = \diagg(\xi, 1,\ldots,1).
$$

It is only left to calculate the order of elements in the periodic
component $g_0H$. As we know, $\forall g\in g_0H:\,\ord(g) =
\ord(g_0)$ equals the least common multiple of $\ord(\Sigma)$ and
$\ord(\Sigma D)$, but since $\Sigma^n = \xi^{n-1}E$, we have
$\ord(\Sigma) = nk$, and
$$(\Sigma D)^n =
\begin{pmatrix}
0 &&&&\xi^{n-1}\\
\xi & 0 &&&\\
& \ddots & \ddots &&\\
&& 1 & 0 &\\
 &&& 1& 0
\end{pmatrix}^{n} =
\xi^nE = E,
$$
so we receive that $\ord(g_0) = nk$. }
\end{examp}

\section{Torus extensions}

Let the connected component of unity of the group $G$ be a torus
$T$. In this case the criterion of the existence of periodic
components may be formulated in a more convenient way. Note that the
automorphism $\varphi_g$ induces the automorphism $A_{\varphi_g}$ of
character lattice of torus $T$:
$$(A_{\varphi_g}\circ \lambda)(t) =
\lambda(\varphi_g(t)).$$

Further, the character lattice of $T$ is a free abelian group with
generators $\chi_s: T \cong (\Bbbk^*)^m\longrightarrow\Bbbk^*,\;
(t_1,\ldots, t_m)\mapsto t_s$. Using the additive notation for the
group of characters we may associate the matrix $A_g\in
GL_n(\mathbb{Z})$ with the automorphism $A_{\varphi_g}$.

\begin{utver} If $G^0 = T$, then the component $gT$ is periodic if
and only if the matrix $A_g$ has no eigenvalue equal to $1$.
\end{utver}

\begin{proof} Let the component $gT$ be periodic and $g^n = e$.
Then by Corollary 1 $\forall h\in T:\,(gh)^n = e$. But evidently
$(gh)^n = g^n\cdot\varphi_g^{n -
1}(h)\cdot\ldots\cdot\varphi_g(h)\cdot h$, and so $\forall h\in T:\,
\varphi_g^{n - 1}(h)\cdot\ldots\cdot\varphi_g(h)\cdot h = e$ or in
terms of the matrix $A_g:$
$$A_g^{n - 1} +\ldots + A_g + E = 0.$$
Consequently the matrix $A_g$ has no eigenvalue equal to 1.

\medskip

Conversely, assume that $g$ is an element of infinite order and that
the matrix $A_g$ has no eigenvalue equal to $1$. Since the number of
connected components of the group $G$ is finite, we have that
$\exists k\in\mathbb{N}:\,g^k\in T$ and hence $\varphi_g^k = id$.
Consequently, $0 = A_g^k - E = (A_g - E)(A_g^{k - 1} + \ldots + A_g
+ E)$. But by our assumption the matrix $A_g - E$ is non-singular,
and so $A_g^{k - 1} + \ldots + A_g + E = 0$ yielding $\forall h\in
T:\,\varphi_g^{k - 1}(h)\cdot\ldots\cdot\varphi_g(h)\cdot h = e$.
But $g^k$ lyes in $T$ and has infinite order and therefore
$\varphi_g^{k - 1}(g^k)\cdot\ldots\cdot\varphi_g(g^k)\cdot g^k =
g^{k^2}\ne e$ which is a contradiction.

\end{proof}

\begin{fiend} {\rm Let $T$ be a torus, $G = T\cup gT\cup\ldots\cup g^{n-1}T$
be its cyclic extension with periodic component $gT$. By Theorem 1
in this case $gT = \{tgt^{-1}\left|t\in T\right.\}$, and hence
$(gT)^n = \{tg^nt^{-1}\left|t\in T\right.\}$ is a conjugacy class of
the element $g^n\in T$. But a torus is commutative yielding $(gT)^n
= \{g^n\}$.

However for general solvable groups it is not true. Return to
Example 2. Since
$$g_0^n =
\begin{pmatrix}
\xi^{n-1}E & 0\\
0 & E
\end{pmatrix}\notin Z_H(U),$$
we have $(g_0H)^n = \{hg_0^nh^{-1}\left|h\in H\right.\}\supsetneqq
\{g_0^n\}$. Hence Example 2 illustrates the difference between
finite extensions of tori and general solvable groups.}
\end{fiend}

The future examples will show that the eigenvalues of the matrix
$A_g$ do not define the order of elements in the component $gT$.
However, we may give some estimates. For this we formulate a
slightly more general proposition.

\medskip

\begin{utver} Let $T$ be a torus,
$G = T\cup gT\cup\ldots\cup g^{m-1}T$ be its cyclic extension of
order $m$ and $k$ be the order of the automorphism $\varphi_g$ where
$\varphi_{g}(t) = g^{-1}tg$. Then there exists an element $g_0\in
gT$ such, that $g_0^{mk}=e$.
\end{utver}

\begin{proof} Let $P_g(t) =
\varphi^{m-1}_g(t)\cdot\ldots\cdot\varphi^{2}_g(t)\cdot\varphi_g(t)\cdot
t$; $Q_g(t) =
\varphi^{k-1}_g(t)\cdot\ldots\cdot\varphi^{2}_g(t)\cdot\varphi_g(t)\cdot
t.$ Then $\forall t\in T: (gt)^m =  g^m\cdot P_g(t)$ и $(gt)^k =
g^m\cdot Q_g(t)$.
 Since $\varphi_g$ is a homomorphism and torus is commutative, $P_g$
 and
 $Q_g$ are homomorphisms of torus $T$. Also $\varphi_g^k = id\Rightarrow P_g(t) =
 (Q_g(t))^r$, where $m = kr.$

 \medskip

\begin{laemma} The following equations are satisfied:

1. $Q_g(Q_g(t)) = (Q_g(t))^k$;

2. $ P_g(P_g(t)) = (P_g(t))^m$.
\end{laemma}

\begin{proof} Using that $\varphi_g$ is a homomorphism and
$\varphi_g^k = id$:
$$Q_g(Q_g(t)) = Q_g(\varphi^{k-1}_g(t)\cdot\ldots\cdot\varphi_g(t)\cdot
t) = $$
$$ = \prod_{i = 0}^{k - 1}\varphi_g^i(\varphi^{k-1}_g(t)\cdot\ldots\cdot
\varphi_g(t)\cdot t) = \prod_{i = 0}^{k -
1}\varphi_g^i(\varphi^{k-1}_g(t))\cdot\ldots
\cdot\varphi_g^i(\varphi_g(t))\cdot \varphi_g^i(t) = $$
$$ = \prod_{i = 0}^{k - 1}\varphi^{k+i-1}_g(t)\cdot\ldots\cdot\varphi_g^{i+(k-i+1)}(t)\cdot\varphi_g^{i+(k-i)}(t)
\cdot\varphi_g^{i+(k-i-1)}(t)\cdot\ldots \cdot \varphi_g^i(t) = $$
$$ = \prod_{i = 0}^{k - 1}\varphi^{i-1}_g(t)\cdot\ldots\cdot\varphi_g(t)\cdot t
\cdot\varphi_g^{k-1}(t)\cdot\ldots \cdot\varphi_g^{i+1}(t)\cdot
\varphi_g^i(t) = (Q_g(t))^k.$$

To prove the second statement note that
$$P_g(P_g(t)) = (Q_g((Q_g(t))^r))^r = (Q_g(Q_g(t)))^{r^2} = $$
$$ = ((Q_g(t))^k)^{r^2} = ((Q_g(t))^r)^{kr} = (P_g(t))^m.$$
\end{proof}

\medskip

\begin{corol} The groups $\Ker{P_g}\cap \TrueIm{P_g}$ and
$\Ker{Q_g}\cap \TrueIm{Q_g}$ are finite.
\end{corol}

\begin{proof} If $s\in \TrueIm{P_g}$, then $\exists\,t\in T: s = P_g(t).$ Further, $s\in
\Ker{P_g} \Rightarrow e = P_g(s) = P_g(P_g(t)) = (P_g(t))^m = s^m.$
But in a torus there are only finitely many elements of order not
greater then $m$. For $Q_g$ the proof is just the same.
\end{proof}

We have $\left.T\right/\Ker{P_g}\cong \TrueIm{P_g}$ and,
consequently, $\dim{\Ker{P_g}} + \dim{\TrueIm{P_g}} = \dim{T}.$
Consider now the homomorphism
$$\gamma : \Ker{P_g}\times\TrueIm{P_g} \longrightarrow T,\quad (t_1,
t_2)\mapsto t_1\cdot t_2.$$

Since $\Ker{\gamma} = \{(t, t^{-1})\}\subset
\Ker{P_g}\times\TrueIm{P_g}$ and is by Corollary 4 a finite subgroup
we have $T = \Ker{P_g}\cdot\TrueIm{P_g}.$

\par
\begin{laemma} $\TrueIm{Q_g} = \TrueIm{P_g}$
\end{laemma}

\begin{proof} If $t\in\TrueIm{P_g}$, then $\exists\, s\in T: t = P_g(s).$
Therefore $Q_g(s^r) = (Q_g(s))^r = P_g(s) = t$, i.e.
$t\in\TrueIm{Q_g}.$
 \medskip
 Conversely, let $t\in\TrueIm{Q_g}$, i.e. $\exists\, s\in T: t =
 Q_g(s)$. We should prove that $\exists\, s_1\in T: t = P_g(s_1).$ But the ground
 field is algebraically closed, therefore $\exists s'\in T: (s')^r = s$,
 implying $P_g(s') = (Q_g(s'))^r = Q_g((s')^r)
  = Q_g(s) = t.$ This means that $t\in\TrueIm{P_g}$.
\end{proof}

Hence $\left.T\right/\Ker{Q_g}\cong \TrueIm{Q_g} = \TrueIm{P_g}$,
 and repeating the previous speculations we get that $T = \Ker{Q_g}\cdot\TrueIm{P_g}$

 \medskip

Now note that the image of the mapping
 $$\alpha : T\longrightarrow T, t\mapsto (gt)^m = g^m\cdot P_g(t)$$
 \noindent equals to the coset $g^m\TrueIm{P_g}\subset T$.
 But $T = \Ker{Q_g}\cdot\TrueIm{P_g}$, and consequently any such
 coset contains an element of $\Ker{Q_g}$, i. e. $\exists\, t_0\in
 T$
such, that $(gt_0)^m = g^m\cdot P_g(t) = s\in\Ker{Q_g}.$ Without
loss of generality we may assume that $g^m\in\Ker{Q_g}$.

 \medskip

\begin{laemma} If $g^m\in\Ker{Q_g}$, then $g^{mk} = e$
\end{laemma}

 \begin{proof} By the definition
 $\forall l=1,\ldots, m-1:\ \varphi_g^l(g^m)=g^{-l}g^mg^l=g^m$.
 But $g^m\in\Ker{Q_g}$ and therefore $e = Q_g(g^m) =
\varphi_g^{k-1}(g^m)\cdot
 \ldots\cdot\varphi_g(g^m)\cdot g^m = (g^m)^k$, yielding $g^{mk} = e$
\end{proof}
Proposition 6 is proved.
\end{proof}

\medskip

\begin{corol} If the component $gT$ is periodic, then the order of
any its element divides $mk$, where $m = \ord(\varphi_g)$ and $k$ is
the order of $gT$ in $G/T$.
\end{corol}

In Example 4 will be shown that there exists a torus $T$ and its
automorphism $\varphi$ such, that any extension $G = T\cup
gT\cup\ldots\cup g^{m-1}T$, in which $\forall t\in T : g^{-1}tg =
\varphi(t)$, satisfies the inequality $\ord(g') < m\ord(\varphi)\
\forall g'\in gT$. But in many situations the estimate proves to be
precise.

\begin{examp} {\rm Consider $G =
T\cup gT\cup\ldots\cup g^{rk-1}T$, where $k$ and $r$ are natural
numbers and $T$ is as follows:
 $$T = \left\{\left.
 \begin{pmatrix}
 t_1 & \ & \ &\ & \ \\
 \ & t_2 & \ &\ & 0\\
 \ & \ & \ddots & \ & \ \\
 \ & \ & \ & t_{k - 1} & \ \\
 0 & \ & \ &\ & t_k\\
 \end{pmatrix}
 \right|t_1t_2\ldots t_{k-1}t_k = 1\right\}\
 \ \subset\ GL_k(\Bbbk),$$
$$ g = \begin{pmatrix}
 0 & \ & \ &\ & \xi\\
 1 & 0 & \ & 0 & \ \\
 \ & 1 & \ddots & \ & \ \\
 \ & \ & \ddots &\ddots & \\\
 0 & \ & \ & 1 & 0\\
 \end{pmatrix},$$
 \noindent where $\xi$ is a primitive root of unity of degree
 $k$. In this case

$$A_g =
\begin{pmatrix}
0 & \ & \ & \ & -1\\
1 & 0 & \ & \ & -1 \\
\ & \ddots & \ddots & \ &\vdots \\
\ & \ & 1 & 0 & -1 \\
\ & \ & \ & 1 & -1 \\
\end{pmatrix}\in GL_{k-1}(\mathbb{Z}).$$

It is easy to see that the eigenvalues of the matrix $A_g$ are the
following: $\eta, \eta^2,\ldots, \eta^{k - 1}$, where $\eta$ is a
primitive root of unity of degree $k$, and therefore the connected
component $gT$ is periodic, and the order of its elements are equal
to $k^2r$. }
\end{examp}

Let now the component $gT$ of cyclic extension $G = T\cup
gT\cup\ldots\cup g^{m-1}T$ be periodic. Since $g^m\in
T^{\varphi_g}$, we have $\ord(g^m)\leqslant \max\{\ord(t)\mid t\in
T^{\varphi_g}\}$. Using this consideration, one can get the
following estimate on the order of $g$.

\begin{utver} For the periodic component $gT$ of cyclic extension
$G = T\cup gT\cup\ldots\cup g^{m-1}T$ the inequality is satisfied:
$\ord(g)\leqslant m|\chi_g(1)|$, where $\chi_g(\lambda)$ is the
characteristic polynomial of the matrix $A_g$.
\end{utver}

\begin{proof} Let $\dim{T} = n$.
Note that fixed points of the automorphism $\varphi_g$ can be found
from the following simultaneous equations:
\begin{equation*}
\begin{cases}
t_1^{a_{11}}t_2^{a_{21}}\ldots t_n^{a_{n1}} = t_1,\\
t_1^{a_{11}}t_2^{a_{22}}\ldots t_n^{a_{n2}} = t_2,\\
\ldots\\
t_1^{a_{1n}}t_2^{a_{2n}}\ldots t_n^{a_{nn}} = t_n,
\end{cases}
\end{equation*}
where $A_g = (a_{ij})$. This means that
\begin{equation*}
\begin{cases}
t_1^{a_{11} - 1}t_2^{a_{21}}\ldots t_n^{a_{n1}} = 1,\\
t_1^{a_{12}}t_2^{a_{22} - 1}\ldots t_n^{a_{n2}} = 1,\\
\ldots\\
t_1^{a_{1n}}t_2^{a_{2n}}\ldots t_n^{a_{nn} - 1} = 1.
\end{cases} \eqno{(1)}
\end{equation*}
By multiplying one equation on another and changing their order we
will be performing the integer elementary transformations upon the
rows of the matrix $B = (a_{ij} - \delta_{ij})$. Such
transformations on one hand do not change the absolute value of its
determinant and on the other hand using them we shall make the
matrix $B$ upper-triangular. System (1) will turn to
\begin{equation*}
\left\{\begin{aligned}
t_1^{a_{11}'}t_2^{a_{12}'}t_3^{a_{13}'}\ldots t_n^{a_{1n}'} &= 1,\\
t_2^{a_{22}'}t_3^{a_{23}'}\ldots t_n^{a_{2n}'} &= 1,\\
\ldots\\
t_n^{a_{nn}'} &= 1.
\end{aligned}\right.
\end{equation*}
So we get that $t_n^{a_{nn}'} = 1,\ t_{n-1}^{a_{n-1,n-1}'} =
t_n^{-a_{n-1,n}'},\ldots, t_1^{a_{11}'} = t_2^{-a_{12}'}\ldots
t_n^{-a_{1n}'}$. Hence while solving the equations we shall at first
extract a root of degree $|a_{nn}'|$ from unity, then a root of
degree $|a_{n-1,n-1}'|$ from a certain root of unity of degree
$|a_{nn}'|$ and so on, but we shall never get more than a root of
degree $|a_{11}'\ldots a_{nn}'|$ of unity. The order of the fixed
point of the automorphism $\varphi_g$ that we have found is equal to
the lesser common multiply of $\ord(t_1),\ldots,\ord(t_n)$, but
$\ord(t_i)$ divides $|a_{ii}'\ldots a_{nn}'|$, and therefore their
lesser common multiply divides $|a_{11}'\ldots a_{nn}'| =
|\det(a_{ij}')| = |\det(B)| = |\det(A - E)| = |\chi_g(1)|$.
\end{proof}

In some situations this estimate is stronger then the estimate of
Corollary 5.

\begin{examp} {\rm For any extension $G = T\cup
gT\cup\ldots\cup g^{m-1}T$ with $\chi_g(\lambda) = \frac{\lambda^r +
1}{\lambda + 1}$, where $r = \dim T + 1$, we have $\chi_g(1) = 1$.
This means that $\ord(g) = m$ yielding $G\cong T\leftthreetimes
\langle g\rangle_m$. }
\end{examp}

\section{The normalizer of a maximal torus of a simple group}

In this section we study the normalizers of maximal torus of
classical simple groups as well as of the exceptional group $G_2$
and partially of $F_4$ and $E_8$, finding their periodic components
and the orders of elements in these components.

First of all we recall a few general facts. Let $T$ be a maximal
torus of a simply connected simple algebraic group $G$, $N_G(T)$ be
its normalizer in $G$, $W = N_G(T)/T$ be its Weyl group, and
$\{\alpha_1,\ldots,\alpha_n\}$ be a system of simple roots. A
product $c = r_1r_1\ldots r_n$, where $r_i$ are the reflections
associated with simple roots, is called a {\it Coxeter element} of
group $W$, and the order $h$ of $c$ is called the {\it Coxeter
number} of $W$. Note that all Coxeter elements are conjugate in $W$,
and that every element of the Weyl group conjugate to a Coxeter
element is also a Coxeter element for a certain system of simple
roots. In particular $c$ is conjugate to $c^{-1}$.

As it is shown in ~\cite[Ch. 3.16]{Reflect}, a Coxeter element has
no eigenvalue equal to $1$. Hence its eigenvalues are as follows:
$\zeta^{m_1},\ldots,\zeta^{m_n}$, where $\zeta$ is a primitive root
of unity of degree $h$ and $m_1\leqslant m_2\leqslant\ldots\leqslant
m_n$. The numbers $m_i$ are called the {\it exponents} of $G$. The
values of the exponents for each simple group can be found, for
example, in~\cite{Vind}. All that was said above immediately proves
the following proposition.

\begin{utver} The component $gT$ of the normalizer $N_G(T)$,
corresponding to a Coxeter element of the Weyl group is periodic.
\end{utver}

As it is shown in~\cite[Prop. 30]{Kill}, a Coxeter element commutes
only with its powers, so the number $N_c$ of Coxeter elements in $W$
equals to $|W|/h$.

However in many simple groups the normalizer of maximal torus has
periodic components corresponding neither to Coxeter elements nor to
their powers. Consider a natural action of the group $W$ in the
rational vector space spanned by the lattice of characters of $T$.
To find the number of all periodic components, we use the following
proposition, see~\cite{Solomon}.

\begin{utver} If $g_k$ is the number of elements of Weyl group $W$,
whose dimension of the space of fixed points is equal to $n-k$, then
$\sum_{k = 0}^{n}{g_kz^k} = \prod_{i = 1}^n{(1 + m_iz)}$.
\end{utver}

By Theorem 1 the number of periodic components in $N_G(T)$ equals to
$g_n = m_1\ldots m_n$. The results are collected in the table:

\medskip

\begin{center}
\begin{tabular}{|c|c|c|c|c|c|}
\hline
 & $SL_n$ & $SO_{2n}$ & $SO_{2n+1}, Sp_{2n}$ & $G_2$ & $F_4$ \\
\hline
$|W|$ & $n!$ & $2^{n-1}n!$ & $2^nn!$ & $12$ & $1152$\\
\hline
$g_n$ & $(n - 1)!$ & $(2n-3)!!(n-1)$ & $(2n-1)!!$ & $5$ & $385$\\
\hline
$N_c$ & $(n-1)!$ & $2^{n-2}(n-2)!n$ & $2^{n-1}(n-1)!$ & $2$ & $96$\\
\hline
\end{tabular}

\smallskip

\begin{tabular}{|c|c|c|c|}
\hline
& $E_6$ & $E_7$ & $E_8$\\
\hline
$|W|$ & $2^{7}\cdot 3^4\cdot 5$ & $2^{10}\cdot 3^4\cdot 5\cdot 7$ & $2^{14}\cdot 3^5\cdot 5^2\cdot 7$\\
\hline
$g_n$ & $2^5\cdot 5\cdot 7\cdot 11$ & $3^2\cdot 5\cdot 7\cdot 11\cdot 13\cdot 17$ & $7\cdot 11\cdot 13\cdot 17\cdot 19\cdot 23\cdot 29$\\
\hline
$N_c$ & $2^{5}\cdot 3^3\cdot 5$ & $2^{9}\cdot 3^2\cdot 5\cdot 7$ & $2^{13}\cdot 3^4\cdot 5\cdot 7$\\
\hline
\end{tabular}
\end{center}

\medskip

Now we find the order of elements in periodic components.

\medskip

{\large 1. The group $SL_n$.}

It is easy to see that Coxeter elements in the Weyl group consisting
of permutations on $n$ elements are precisely the cycles of length
$n$. Their number is equal to $(n-1)!$ and it is the number of
periodic components in $N(T)$.

To calculate the order of elements in the component corresponding to
a cycle of length $n$, consider the monomial $(0,1, -1)$-matrix $A$,
that lyes in it. Here
$$\ord(A) =
\begin{cases}
n&\text{for odd $n$};\\
2n&\text{for even $n$}.
\end{cases}$$

\medskip

{\large 2. The group $SO_{2n}$.}

Consider the following matrix representation: $SO_{2n}\cong \{g\in
GL_{2n}\mid g^T\Omega g = \Omega\}$, where $\Omega$ is a monomial
matrix with unities on the secondary diagonal.

The normalizer of a maximal torus $T$ is
$$\left\{A = (a_{i,j})\in SL_{2n}\left|\begin{matrix}
\ A \mbox{ is a monomial matrix};\\ a_{2n - i + 1, 2n - j + 1} =
a_{i, j}^{-1}\ \mbox{если}\ a_{i,j}\ne 0
\end{matrix}
\right.\right\}.$$

Consider a component $gT$ of the normalizer. It defines a
homomorphism $\varphi_g : t\mapsto gtg^{-1}$, and the latter has a
corresponding matrix $A_g$ of the mapping induced in the character
lattice. It is easy to see that it is a monomial $(0,1,-1)$-matrix.
We call a monomial matrix $A\in GL_n(\mathbb{Z})$ corresponding to a
permutation $\tau\in S_n$ (we denote it by $A = A_{\tau}$), if
$\forall j = 1,\ldots,n\,\exists \alpha_j\in \Bbbk^*:\,Ae_j =
\alpha_je_{\tau(j)}$. Then the matrix $A_g$ corresponds to a certain
permutation $\sigma\in S_n$.

Let $\sigma = \sigma_1\ldots\sigma_k$ be the decomposition in a
product of independent cycles (and of length 1 as well). Then $A$ is
conjugate to the matrix
$$A' = \begin{pmatrix}
S_{\sigma_1} &&&\\
& S_{\sigma_2} &&\\
&& \ddots &\\
&&& S_{\sigma_k}
\end{pmatrix},$$
where $S_{\sigma_j}$ is a monomial matrix corresponding to
$\sigma_j$. But if $\tau\in S_m$ is a cycle of length $m$, then
$(-1)^{\tau} = (-1)^{m - 1}$ and the characteristic polynomial of
$S_{\tau}$ is equal to
$$\chi_{S_{\tau}}(\lambda) = \lambda^m + (-1)^{\tau}(-1)^kP
= \lambda^m - P,$$ where $P$ is the product of non-zero elements of
$S$. Hence the characteristic polynomial of the matrix $A$ equals to
$$\chi_A = \chi_{S_{\sigma_1}}\cdot\ldots\cdot\chi_{S_{\sigma_k}}
= (\lambda^{s_1} - P_1)\ldots(\lambda^{s_k} - P_k),$$ where $s_j$ is
the order of $\sigma_j$ and $P_j$ is the product of non-zero
elements of $S_{\sigma_j}$.

It is clear now that the matrix $A$ has no eigenvalue equal to 1 if
and only if all $P_j\ne 1$, and since non-zero elements of $A$ may
be only equal to $\pm 1$ (and of $A'$ as well, because it has been
made from $A$ by permuting rows and columns), we get that the matrix
$A$ has no eigenvalue equal to 1 if and only if in each block
$S_{\sigma_j}$ of the matrix $A'$ the number of elements -1 is odd.

Finally, periodic components correspond to the following elements of
the Weyl group: $(\sigma, \theta)$, where $\sigma\in S_n$, $\theta =
(\varepsilon_1,\ldots,\varepsilon_n)\in\{\pm 1\}^n$, and if $\sigma
= (i_{1,1},\ldots,i_{1,s_1})$ $\ldots(i_{k,1},\ldots,i_{k,s_k})$ is
the decomposition in a product of independent cycles (of length 1 as
well), then
$$\forall j = 1,\ldots,k:\,\prod^{s_j}_{q = 1}\varepsilon_{i_{j,q}} = -1.\eqno{(2)}$$
(Note that the number of minus unities among the components of the
vector $\theta$ is even).

Now we calculate the order of elements in the corresponding
component. For that consider the monomial $(0,1)$-matrix $g$ lying
in it. Direct calculations show that the order of $g$, and with it
the order of every element of the component equals to double least
common multiple of $s_1,\ldots,s_k$.

\pagebreak

{\large 3. The group $SO_{2n+1}$.}

The results are the same as in the paragraph 2. The only difference
is that the number of minus unities among the components of the
vector $\theta$ need not be even.

\medskip

{\large 4. The group $Sp_{2n}$.} \nopagebreak

Consider the following matrix representation: $Sp_{2n}\cong \{g\in
GL_{2n}\mid g^T\Lambda g = \Lambda\}$, where
$$\Lambda =
\begin{pmatrix}
0 & I\\
-I & 0
\end{pmatrix}$$
(here $I$ is a monomial matrix with unities on the secondary
diagonal).

The normalizer of maximal torus $T$ is
\begin{equation*}
\left\{A = (a_{i,j})\in GL_{2n}\left|\
\begin{split}
&A \mbox{ is a monomial matrix};\\
& a_{2n - i + 1, 2n - j + 1} =
a_{i, j}^{-1},\\
&\ \ \mbox{if}\ a_{i,j}\ne
0\ \mbox{and}\ i,j\leqslant n\ \mbox{or}\ i,j > n;\\
&a_{2n - i + 1, 2n - j + 1} = - a_{i, j}^{-1},\\
&\ \ \mbox{if}\ a_{i,j}\ne 0\ \mbox{and}\ i\leqslant n < j\
\mbox{or}\ j\leqslant n < i
\end{split}
\right.\right\}.
\end{equation*}

It is easy to see that the Weyl group here is the same as for
$SO_{2n + 1}$, so we can use the results of the paragraph 2.
Periodic components correspond to the following elements of the Weyl
group: $(\sigma, \theta)$, where $\sigma\in S_n$, $\theta =
(\varepsilon_1,\ldots,\varepsilon_n)\in\{\pm 1\}^n$, so that if
$\sigma = (i_{1,1},\ldots,i_{1,s_1})$
$\ldots(i_{k,1},\ldots,i_{k,s_k})$ is the decomposition in a product
of independent cycles (of length 1 as well), then the equations (2)
are satisfied.

Direct calculations show that the order of the monomial
$(0,1,-1)$-matrix lying in $gT$ equals to the least common multiple
of $s_1,\ldots,s_k$ multiplied by 4.

\medskip

{\large 5. The group $G_2$.}

The corresponding Weyl group consists of 12 elements: the identity,
six reflections and five rotations through angles multiple to
$\frac{\pi}{6}$. The rotations are the powers of Coxeter element $c$
and have no eigenvalue equal to 1. Hence periodic components
correspond to the elements $c,\ldots,c^5$.

To find the orders of elements in these components consider the
characteristic polynomial $\chi_c(\lambda)$ of $A_c$. It is
quadratic, its coefficients are integers and the roots are equal to
$\zeta$ and $\zeta^5$, where $\zeta$ is a primitive root of unity of
degree $6$. One can easily understand that $\chi_c(\lambda) =
\lambda^2 - \lambda + 1$. Proposition 7 yields then that for every
element $g$ of the component corresponding to $c$, the inequality
$\ord(g)\leqslant 6|\chi_c(1)| = 6$ is satisfied, implying $\ord(g)
= 6$. The orders of elements in the components corresponding to
$c^2,\ldots,c^5$ are equal to 3, 2, 3 and 6. 

\medskip

{\large 6. The group $F_4$.}

Here we find the order of elements in the components corresponding
to Coxeter element and its powers.

Note that the exponents of $F_4$ are equal to 1, 5, 7 and 11 and
hence are coprime with the Coxeter number $h = 12$. Therefore all
the degrees of Coxeter elements correspond to periodic components.

Let the component $gT$ correspond to Coxeter element $c$. By
Corollary 5 we have $\ord(g) = \ord(c)a = 12a$, where $a|12$. But
then $\ord(g^6) = 2a$, and by Corollary 5 we get $a|2$. Considering
$\ord(g^4)$, we get $a|3$, which means that $a = 1$. Finally the
orders of elements in the components corresponding to
$c,c^2,c^3,c^4,c^5,c^6$, are equal to 12, 6, 4, 3, 12, 2. It is easy
to see that in the normalizer of maximal torus of the group $F_4$
there are periodic components corresponding neither to Coxeter
elements nor to their powers.

Note that similar speculations may be performed for Coxeter elements
of the group $E_8$. But as for the elements of periodic components
of the normalizer of maximal torus of the groups $E_6$ and $E_7$,
and of periodic components of the normalizer of maximal torus of
$F_4$ and $E_8$, not corresponding to powers of Coxeter elements,
calculation of their order requires further investigation.

\end{document}